\input amstex
\magnification=\magstep1 \baselineskip=13pt
\documentstyle{amsppt}
\vsize=8.7truein \CenteredTagsOnSplits \NoRunningHeads
\def\today{\ifcase\month\or
  January\or February\or March\or April\or May\or June\or
  July\or August\or September\or October\or November\or December\fi
  \space\number\day, \number\year}

\def\EE{{\bold E\thinspace }}
\def\conv{\operatorname{conv}}

\def\vl{\operatorname{vol}}
\def\Mat{\operatorname{Mat}}

\def\co{\operatorname{co}}
\topmatter
\title The Computational Complexity of Convex Bodies \endtitle
\author Alexander Barvinok  and Ellen Veomett \endauthor
\address Department of Mathematics, University of Michigan, Ann Arbor,
MI 48109-1043 \endaddress \email barvinok$\@$umich.edu, eveomett$\@$umich.edu
\endemail
\date October  2006 \enddate
\thanks This research was partially supported by NSF Grant DMS 0400617.
\endthanks
\abstract We discuss how well a given convex body $B$ in a real $d$-dimensional vector space $V$
 can be approximated by a set 
$X$ for which the membership question: ``given an $x \in V$, does $x$ belong
to $X$?'' can be answered efficiently (in time polynomial in $d$).
We discuss approximations of a convex 
body by an ellipsoid, by an algebraic hypersurface,  by a projection of a polytope with a controlled number of facets, and by a section of the cone of positive semidefinite quadratic forms. We illustrate some of 
the results on the Traveling Salesman Polytope, an example of a complicated convex body studied in  
combinatorial optimization. 
\endabstract
\endtopmatter
\document

\head 1. Introduction \endhead

Let $V$ be a finite-dimensional real vector space. A set $B \subset V$ is called {\it convex}
if for every two points $x, y \in B$ the interval $[x, y]=\bigl\{ \alpha x + (1-\alpha)y: \ 0 \leq \alpha \leq 1
\bigr\}$ also lies in $B$. A convex set $B \subset V$ is a {\it convex body} if $B$ is compact 
with a non-empty interior. 

Various notions of complexity for convex bodies have been discussed from a variety of points 
of view, see, for example, \cite{Sz06} for a recent survey.   In this paper, we suggest the following 
general approach to the {\it computational} complexity of convex bodies.

Let $B \subset V$ be a $d$-dimensional convex body. We would like to design an efficient 
algorithm for the following {\it membership problem}

\specialhead (1.1) The membership problem for $B$ \endspecialhead
\medskip
{\bf Input:} A point $x \in V$.
\medskip
{\bf Output:} ``Yes'' if $x \in B$ and ``No'' if $x \notin B$.
\bigskip

The algorithm in Problem 1.1 is, of course, determined by the convex body $B$. It is important 
to note that we do not count the time and resources spent on creating such an algorithm 
towards the computational complexity of $B$. Once the algorithm is created, only the time 
needed to answer the question whether a given point $x \in V$ belongs to $B$ counts as the 
complexity of $B$.

At this point, it is not important to us what model of computation we use. For example, 
we can identify $V={\Bbb R}^d$ and assume that a point $x=\left(\xi_1, \ldots, \xi_d\right)$ 
is given by its real coordinates. The algorithm employs arithmetic operations and its complexity 
is the number of operations used (the real model). Alternatively, we may assume that the point 
$x$ has rational 
coordinates and count bit operations instead (the bit model). 
We are interested in how the complexity of the algorithm 
depends on the dimension $d$.

The algorithmic theory of convexity, see \cite{G+93}, implies that as soon as there is a 
polynomial time algorithm in the Membership Problem 1.1 for a convex body $B$ (in the bit model), there is a polynomial time algorithm in the problem of optimizing a given linear function on $B$
(again, in the bit model).
Hard optimization problems supply a variety of convex bodies for which the membership problem
is likely to be hard. We concentrate on one model example, the (symmetric) {\it Traveling Salesman 
Polytope}.

\example{(1.2) Example: the Traveling Salesman Polytope} Let us fix an integer $n \geq 4$.
In the space $\Mat_n$ of $n \times n$ real matrices, we define $TSP_n$ as the convex 
hull of the adjacency matrices of the Hamiltonian cycles in a complete graph on $n$ vertices.
That is, for any Hamiltonian cycle $C$ (a cycle that visits every vertex exactly once) of a 
complete undirected graph on the vertices $\{1, \ldots, n\}$, we introduce the $n \times n$ 
matrix $x=x(C)$ where 
$$x_{ij}=\cases 1&\text{if\ } ij \ \text{is an edge of} \ C \\
0& \text{otherwise} \endcases$$
and take the convex hull of all such matrices $x(C)$. One can observe that $TSP_n$ has 
$(n-1)!/2$ vertices and that $\dim TSP_n=(n^2-3n)/2$, see Chapter 58 of \cite{Sc03}.
We define the ambient space $V_n$ as the affine hull of $TSP_n$ with the origin at the center 
$a=\left(a_{ij}\right)$, where 
$$a_{ij}=\cases {2\over n} &\text{if\ } i \ne j \\ 0 &\text{if \ } i=j \endcases$$
The space $V_n$ consists of the $n \times n$ symmetric matrices with zeros on the diagonal and 
the row and column sums equal to 2.
\endexample

The facial structure of $TSP_n$ was studied in many papers but still remains a mystery, see
Chapter 58 of \cite{Sc03} for a survey.
In some sense, unless the complexity hierarchy collapses, the facets of $TSP_n$ 
{\it cannot} be described. More precisely, there can be no algorithm of a polynomial in $n$
complexity, which, given an inequality in $V_n$ decides whether it defines a facet of $TSP_n$, see
Chapter 4 of \cite{Sc03}.

Since the Membership Problem 1.1 can be quite hard, we also consider the computational 
complexity of the {\it approximation problem}: we want to find a set $X \subset V$ approximating 
$B$ such that the membership problem for $X$ can be solved within a given complexity,
preferably in time polynomial in $\dim V$.

\specialhead (1.3) Measuring the quality of approximation \endspecialhead

In many important cases, the convex body $B$ is {\it symmetric}, that is $B=-B$ 
and we measure how well $X$ approximates $B$ by a number $\alpha \geq 1$ such that 
$$X \subset B \subset \alpha X.$$ Often, $B$ is not symmetric but has a natural center nevertheless,
as in the case of the Traveling Salesman Polytope. In such a case, taking the center of $B$ for the origin,
we measure the quality of approximation in the same way. In many cases, but not in all, 
the set $X$ is convex as well.
\bigskip
This paper is meant to be a survey although it contains some new results, mostly in 
Section 5 and in Theorems 4.6 and 6.2. The paper is structured as follows.
\medskip
In Section 2, we discuss how well a symmetric convex body can be approximated by an 
ellipsoid, or, more generally, by an algebraic hypersurface of a given degree. In particular, for 
any $\epsilon>0$ and any symmetric convex body $B$ we obtain a set $X$ such that 
$X \subset B \subset \alpha X$ for $\alpha=\epsilon \sqrt{\dim B}$ and the membership problem for $X$ can be solved in time polynomial in $\dim B$.
\smallskip
In Section 3, we review results on approximation of convex bodies by polytopes with 
a controlled number of vertices or facets. We show that to approximate a $d$-dimensional symmetric convex body within a constant factor, it is sufficient to take a polytope with $e^{O(d)}$ vertices.
The exponential in $d$ number of vertices is also necessary in the worst case, as shows the example
of the Euclidean ball.
\smallskip
In Section 4, we discuss approximations of convex bodies by {\it projections} of polytopes with a
controlled number of facets. Note that if $X$ is a projection of a polytope with $N$ facets  then the membership problem for $X$ is a linear programming problem that can be solved in time polynomial in $N$ 
(in the bit model), cf. \cite{G+93}.  On the other hand,
the number of facets of a polytope can grow exponentially if 
the operation of projection is applied, so potentially we get more flexibility for approximation. We describe several interesting phenomena
here. First, we show that to obtain a good approximation of a {\it symmetric} convex body by a projection of 
a polytope with not too many facets, we may have to require the polytope to be {\it non-symmetric}.
This is the case, for example, when the body is the cross-polytope (octahedron).
Then we  discuss an amazing approximation, constructed by A. Ben-Tal and  A. Nemirovski
\cite{BN01} of the Euclidean ball in ${\Bbb R}^d$ within a factor of $(1+\epsilon)$ by the projection of a polytope with only $O\left(d \ln \epsilon^{-1}\right)$ facets. 
Finally, we discuss a general way to approximate an arbitrary convex body by a projection 
of a polytope with a controlled number of facets. As an illustration, we state the 
result of the second author that for any $\epsilon>0$ the Traveling Salesman Polytope $TSP_n$ can be approximated within a factor of $\epsilon n$ by the projection 
of a polytope with $n^{O(1/\epsilon)}$ facets.

\smallskip
In Section 5, we present a construction which, for any $d$-dimensional convex body $B \subset V$  ``almost approximates''  its polar $B^{\circ} \subset V^{\ast}$ by a section of a polytope with not more than $e^{O(\sqrt{d} \ln d)}$ vertices. We consider $V^{\ast}$ embedded into the space $C(B)$ of continuous functions on $B$, construct a polytope $R \subset  C\left(B\right)$ with
$e^{O(\sqrt{d} \ln d)}$ vertices and show that we obtain a good approximation of $B^{\circ}$ if 
we slightly ``bend'' $V^{\ast} \subset C(B)$
and intersect in with $R$. 

\smallskip
In Section 6, we discuss approximations by a section of the cone of positive semidefinite 
quadratic forms. Such approximations also lead to efficient algorithms for the membership problem.
We review by now classic results on approximations of the cut polytope and its variations.
We also present a general construction and illustrate it by the result of the second author 
that the polar $TSP_n^{\circ}$ of the Traveling Salesman Polytope
 can be approximated within a factor of 
$\epsilon n$ by a section of the cone of positive semidefinite quadratic forms in $n^{O(1/\epsilon)}$ variables. 

\head 2. Approximation by algebraic hypersurfaces \endhead

An {\it ellipsoid} $E \subset V$ is a set defined as
$$E=\Bigl\{v \in V: \quad q(v-v_0) \leq 1 \Bigr\},$$
where $q: V \longrightarrow {\Bbb R}$ is a positive definite quadratic form and $v_0 \in V$ is 
a particular point called the center of the ellipsoid. As is known, see for example, \cite{Ba97},
for any convex body $B \subset V$ there is a unique ellipsoid $E \subset B$ that has the largest 
volume among all ellipsoids contained in $B$. It is called the {\it John ellipsoid} of $B$ and it 
satisfies 
$$E \subset B \subset dE,$$
where $d=\dim V$ and the origin of $V$ is moved to the center of $E$.

If $B$ is symmetric then the John ellipsoid $E$ is necessarily centered at the origin and 
satisfies 
$$E \subset B \subset \sqrt{d} E.$$
For the Traveling Salesman Polytope $TSP_n$ the ellipsoid $E$ is centered at the center of 
$TSP_n$, touches the facets of $TSP_n$ defined by the equations $x_{ij}=0$ and satisfies 
$$E \subset TSP_n \subset {(n-3) \sqrt{n} \over 2} E \quad \text{for} \quad n \geq 5,$$
see \cite{BB05}.

Since an ellipsoid $E$ is defined by one quadratic inequality, the membership problem for 
$E$ can be solved in $O\left(d^2\right)$ time for $d=\dim V$. One can ask whether better 
bounds can be obtained by using polynomial inequalities of higher order. 

For symmetric convex bodies such a result was obtained in \cite{Ba03}. We state some general
prerequisites first.

\subhead (2.1) Symmetric convex bodies and norms. Polarity \endsubhead
With a symmetric convex body $B \subset V$ one naturally associates a norm $\| \cdot \|_B$ on $V$
defined by
$$\| v\|_B=\inf\Bigl\{ \lambda >0: \quad v \in \lambda B \Bigr\}.$$
Approximating $B$ by an efficiently computable set $X$ is equivalent to approximating $\|\cdot\|_B$ by 
an efficiently computable function $f$. 

Let $V^{\ast}$ be the space of all linear functions $\ell: V \longrightarrow {\Bbb R}$.
We  recall that the {\it polar} $B^{\circ} \subset V^{\ast}$ of a convex body $B \subset V$
is defined by 
$$B^{\circ}=\Bigl\{ \ell \in V^{\ast}: \quad \ell(v) \leq 1 \quad \text{for all} \quad v \in B \Bigr\}.$$
The standard duality result states that $\left(B^{\circ}\right)^{\circ}=B$ if $B$ is a closed 
convex set containing the 
origin.
\bigskip
The following result was proved in \cite{Ba03}.

\proclaim{(2.2) Theorem} For any symmetric convex body $B \subset V$ and any integer 
$k \geq 1$ there exists a homogeneous polynomial $p: V \longrightarrow {\Bbb R}$ such that 
$p$ is the sum of squares of homogeneous polynomials of degree $k$ and
$$p^{1/2k}(v) \leq \|v\|_B \leq \alpha(d,k)p^{1/2k}(v) \quad \text{for all} \quad 
v \in B,$$
where 
$$\alpha(d,k)={d+k-1 \choose k}^{1/2k} \quad \text{and} \quad d=\dim B.$$
\endproclaim

In other words, the set 
$$X=\Bigl\{v \in V: \quad p(v) \leq 1 \Bigr\}$$
approximates $B$ within a factor of $\alpha(d,k)$.
Note that the membership problem for $X$ can be solved in $O\left(d^{2k}\right)$ time.
If $d \gg k \gg 1$, applying Stirling's formula, we get 
$$\alpha(d,k)\approx  {\sqrt{d} \over \left(k!\right)^{1/2k}} \approx \sqrt{de \over k}.$$
It follows then that for any fixed $\epsilon >0$ we can choose $k=k(\epsilon)$ large enough so that 
$X$ approximates $B$ within a factor of $\epsilon \sqrt{\dim B}$ and the membership problem for 
$X$ can be solved in time  polynomial in $\dim B$. It is not clear whether the set $X$ can always 
be chosen convex, although this is the case in many situations.

\demo{Sketch of proof of Theorem 2.2} Let $V^{\ast}$ be the dual space of all linear functionals 
$\ell: V \longrightarrow {\Bbb R}$ and let 
$$B^{\circ}=\left\{ \ell \in V^{\ast}: \quad \ell(v) \leq 1 \quad \text{for all} \quad v \in B \right\}$$
be the polar of $B$. By the standard duality argument, we can write
$$\|v\|_B=\max\Bigl\{\ell(v): \quad \ell \in B^{\circ} \Bigr\}.$$

 Let $V^{\otimes k}$ be the tensor product of $k$ copies of $V$, so
$\left(V^{\otimes k}\right)^{\ast}=\left(V^{\ast}\right)^{\otimes k}$. Then 
$$\|v\|_B^k=\max\Bigl\{\ell^{\otimes k}\left(v^{\otimes k} \right): \quad \ell \in B^{\circ} \Bigr\}.$$
Now, let us define a symmetric convex set $C \subset \left(V^{\ast}\right)^{\otimes k}$ as the 
convex hull
$$C=\conv\Bigl(\ell^{\otimes k}, -\ell^{\otimes k}: \quad \ell \in B^{\circ} \Bigr).$$
Then 
$$\|v\|_B^k=\max\Bigl\{L\left(v^{\otimes k}\right): \quad L \in C \Bigr\}.$$
Finally, we approximate $C$ by an ellipsoid $E \subset C$. The crucial consideration is 
that since $C$ lies in the symmetric part of $\left(V^{\ast}\right)^{\otimes k}$, we have 
$$\dim C \leq {\dim V + k-1 \choose k}$$
and hence $E$ approximates $C$ within a factor of $\alpha^k(d,k)$.

It remains to notice that the formula 
$$\sqrt{p(v)}=\max\Bigl\{L\left(v^{\otimes k}\right): \quad L \in E \Bigr\}$$
indeed defines a polynomial $p: V \longrightarrow {\Bbb R}$ of degree $2k$ which is a sum of 
squares of polynomials of degree $k$.
{\hfill \hfill \hfill} \qed
\enddemo

There are several open questions related to Theorem 2.2.

\subhead \nofrills{(2.3) The best approximation?} \endsubhead

Suppose that we want to approximate an arbitrary norm $\|v\|_B$ by an expression 
$p^{1/2k}(v)$, where $p: V \longrightarrow {\Bbb R}$ is a homogeneous polynomial 
of degree $2k$. Is the coefficient $\alpha(d,k)$ of Theorem 2.2
best that we can get? This is so in the case of $k=1$, since the $d$-dimensional cube 
$$I_d=\Bigl\{\left(\xi_1, \ldots, \xi_d\right): \quad |\xi_i| \leq 1 \quad \text{for} \quad 
i=1, \ldots, d \Bigr\}$$
or the $d$-dimensional cross-polytope (octahedron)
$$O_d=\Bigl\{\left(\xi_1, \ldots, \xi_d\right): \quad \sum_{i=1}^d |\xi_i| \leq 1 \Bigr\}$$
cannot be approximated by an ellipsoid better than within a
factor of $\sqrt{d}$, cf. \cite{Ba97}.

Similarly, in the non-symmetric situation, the $d$-dimensional simplex
$$\Delta_d =\Bigl\{ \left(\xi_1, \ldots, \xi_{d+1} \right): \quad \sum_{i=1}^{d+1} \xi_i=1
\quad \text{and} \quad \xi_i \geq 0 \quad \text{for} \quad i=1, \ldots, d+1 \Bigr\}$$
cannot be approximated by an ellipsoid better than within a factor of $d$.

For $k >1$ it is not clear whether the bound of Theorem 2.2 is optimal. One can show, however,
that the octahedron $O_d$ cannot be approximated by a hypersurface of degree 4 or 6 better than 
within a certain factor $c\sqrt{d}$ for some absolute constant $c>0$ \cite{Ba03}. The 
idea of the proof is as follows. Once a convex 
body $B$ and the degree $2k$ are fixed, there exists the best approximating polynomial $p$.
Furthermore, the polynomial $p$ can be chosen to be invariant under the group of the symmetries
of $B$. In the case of $O_d$, it follows that one can choose $p$ to be a symmetric 
polynomial in $\xi_1^2, \ldots, \xi_d^2$, from which the estimate can be deduced.

It is interesting to note that while the cube $I_d$ and the octahedron $O_d$ demonstrate 
similar behavior with respect to ellipsoidal approximations, their behavior with respect to 
higher degree approximations very much differ. Indeed, the cube $I_d$ can be approximated
by a hypersurface of degree $2k$ within a factor of $d^{1/2k}$. For that, one can choose
$$p(x)=\sum_{i=1}^d \xi_i^{2k} \quad \text{for} \quad x=\left(\xi_1, \ldots, \xi_d\right).$$

\subhead \nofrills{(2.4) Convex approximation?} \endsubhead

Although in many special cases the function $p^{1/2k}(\cdot)$ constructed in Theorem 2.2 
turns out to be a norm, there is no reason to believe that it is always a norm (but there are 
no explicit counterexamples either). Equivalently, there is no apparent reason why the 
set 
$$X=\Bigl\{v \in V: \quad p(v) \leq 1 \Bigr\}$$
should be convex. It is an interesting question whether one can choose $p$ to be a norm 
and still get the same order of approximation. One natural candidate would be 
$$p(v)=\int_{B^{\circ}} \ell^{2k}(v) \ d \mu(\ell), \tag2.4.1$$
where $\mu$ is a certain probability measure on the polar $B^{\circ}$. One can show that for the 
normalized Lebesgue measure $\mu$ 
 the bounds generally are much weaker than 
those of Theorem 2.2. A natural candidate is the {\it exterior angle} measure defined as follows.

Let ${\Bbb S}^{d-1} \subset {\Bbb R}^d$ be the unit sphere in Euclidean space ${\Bbb R}^d$
endowed with the scalar product $\langle \rangle$  and let 
$\nu$ be the Haar probability measure on ${\Bbb S}^{d-1}$. Let $K \subset {\Bbb R}^d$ be 
a convex body. For a closed subset $X \subset K$, let us define 
$$\mu_K(X)=\nu\Bigl\{c \in {\Bbb S}^{d-1}: \quad \max_{x \in X}
\langle c, x \rangle =\max_{x \in K} \langle c,x \rangle \Bigr\}.$$  
In words: the measure of a closed subset $X$ is the proportion of linear functions of length 1 
that attain their maximum on $K$ at a point $x \in X$. Clearly, the probability measure $\mu_K$ 
depends on the convex body $K$ as well as on the Euclidean structure in ${\Bbb R}^d$.

Although no bounds have been proved, it seems plausible that by choosing a scalar product 
in $V$ in which the John ellipsoid $E$ of $B$ is the unit ball and letting $\mu=\mu_{B^{\circ}}$
in (2.4.1), one obtains a norm $p^{1/2k}$ which gives a similar order approximation as 
in Theorem 2.1. 

\head 3. Approximations by polytopes \endhead

Any convex body $B \subset V$ can be arbitrarily well approximated by the convex hull of a sufficiently 
large finite subset $X \subset B$. Consequently, the membership problem 
for $B$ is replaced by the membership problem for the polytope $P=\conv(X)$ whose 
complexity (in the bit model) is polynomial in the cardinality $|X|$ of $X$. In the dual setting, any convex body can be arbitrarily 
well approximated by 
the intersection of a sufficiently large set $X$ of halfspaces.  Then the complexity 
of the membership problem is linear in $|X|$ in both the real and the bit models.   How large should $X \subset B$ be so that 
$\conv(X)$ approximates $B$ reasonably well? 

The following result is well-known, see for example, Lemma 4.10 of \cite{Pi89}.

\proclaim{(3.1) Lemma}  Let $B \subset V$ be a symmetric convex body. Then, for any 
$1>\epsilon>0$ there is a subset $X \subset B$ such that 
$$|X| \leq \left(1 + {2 \over \epsilon} \right)^d \quad \text{for} \quad d=\dim V$$
and 
$$P \subset B \subset {1 \over 1-\epsilon}P \quad \text{for} \quad P=\conv(X).$$
\endproclaim
\demo{Proof} Let $\| \cdot\|_B$ be the norm associated with $B$, see Section 2.1.
Let $X \subset B$ be a maximal (under inclusion) subset such that 
$$\|x_1-x_2\|_B > \epsilon \quad \text{for all distinct} \quad x_1, x_2 \in X.$$ 
Since $X$ is maximal, for every $x \in B$ there is a point $y \in X$ and a point $z \in \epsilon B$
such that $x-y=z$. In other words,
$$B \subset X + \epsilon B \subset P + \epsilon B.$$
Iterating, we get 
$$B \subset P + \epsilon P + \epsilon^2 P + \ldots + \epsilon^k P + \epsilon^{k+1}B$$
and taking the limit we conclude that 
$$B  \subset {1 \over 1-\epsilon} P.$$
Next, we estimate the number of points in $X$. We notice that 
$$\left(x_1 + {\epsilon \over 2} B \right) \cap \left(x_2 + {\epsilon \over 2} B \right) =\emptyset
\quad \text{for distinct} \quad x_1, x_2 \in X.$$
Furthermore,
$$\bigcup_{x \in X} \left(x + {\epsilon \over 2} B \right) \subset \left(1+ {\epsilon \over 2} \right) B.$$
Hence
$$\vl \left(1+{\epsilon \over 2} \right) B \geq |X| \vl \left({\epsilon \over 2} B \right),$$
from which
$$|X| \leq \left({2+\epsilon \over \epsilon} \right)^d$$
as required.
{\hfill \hfill \hfill} \qed
\enddemo
Although in many cases the bound can be slightly sharpened, the bottom line is that the bound 
on the number $|X|$ of points is exponential in the dimension. In fact, simple volume estimates 
show that even if $B \subset {\Bbb R}^d$ is the Euclidean ball, the number of points needed
to approximate $B$ within a constant factor is exponentially large in $d$.

The following result is proved, for example, in \cite{Ba97}.
\proclaim{(3.2) Theorem} Let $B \subset {\Bbb R}^d$ be the unit ball and let 
$X \subset B$ be a set such that for $P=\conv(X)$ we have 
$$P \subset B \subset \alpha P.$$
Then 
$$|X| \geq \exp\left\{ d \over 2 \alpha^2 \right\}.$$
Similarly, if $P$ is the intersection of a set $X$ of halfspaces such that 
$$P \subset B \subset \alpha P$$ then 
$$|X| \geq \exp\left\{ d \over 2 \alpha^2 \right\}.$$ 
\endproclaim

\head 4. Approximation by projections \endhead

Let $B \subset V$ be a convex body. Suppose that we manage to construct a polytope 
$P \subset W$, where $W$ is some other vector space, and a linear transformation 
$T: W \longrightarrow V$ such that the image $Q=T(P)$ approximates $B$ reasonably well.
Given an $x \in V$, testing whether $x \in Q$ reduces to testing whether the affine subspace
$T^{-1}(x)$ has a non-empty intersection with $P$, which is a linear programming problem.
In particular, if $P$ is defined by $N$ linear inequalities, the complexity of the membership 
problem for $Q$ is bounded by a polynomial in $N$ (in the bit model). On the other hand, the 
number of facets of $Q$ can be exponentially large in $N$, which shows that, in principle, 
$Q$ may provide a fairly good approximation even to rather ``non-polytopal'' convex bodies 
$B$.

In the dual setting, we are interested in approximating the body $B$ by a section of a polytope.
Namely, we want to construct a vector space $W \supset V$ and a polytope $P \subset W$ such that 
the intersection $Q=P \cap V$ approximates $B$ reasonably well. If $P$ is defined as the convex hull
of $N$ vertices, the membership problem for $P$, being a linear programming 
problem, has the complexity that is polynomial in $N$ (in the 
bit model). On the other hand, the number of vertices of $Q=P \cap V$ can be exponentially 
large in $N$, which, in principle, can lead to good approximations of $B$.

One can observe that a polytope $Q$ is a projection of a polytope with at most $N$ facets if 
and only if $Q$ is a section of a polytope with at most $N$ vertices, so the two approaches are essentially identical.

\subhead (4.1) Good sections are good projections and vice versa \endsubhead

Let $V \subset W$ be a pair of spaces, let $P \subset W$ be a polytope with $N$ vertices 
and let $Q=P \cap V$. Since $P$ has $N$ vertices, there exists a simplex 
$\Delta \subset {\Bbb R}^N$, see Section 2.3, and 
a linear transformation $T: {\Bbb R}^N \longrightarrow W$ such that $P=T(\Delta)$. 
Then $U=T^{-1}(V) \subset {\Bbb R}^N$ is a subspace, the polytope $P'=\Delta \cap U$ 
has at most $N$ facets and $Q=T(P')$. In other words, a polytope that is a section of 
a polytope with at most $N$ vertices can be represented as a projection of a polytope 
with at most $N$ facets.

Vice versa, let $P \subset W$ be a polytope with $N$ facets, let $T: W \longrightarrow V$
be a linear transformation and let $Q=T(P)$. Fixing a scalar product in $W$, we identify $V$ 
with a subspace $V \subset W$ and $T$ with the orthogonal projection of $W$ onto $V$.
Assuming that $P$ contains the origin in its interior, let us consider the polars $P^{\circ}$ and 
$Q^{\circ}$. One can see that $Q^{\circ}=P^{\circ} \cap V$. Since $P^{\circ}$ is a polytope with 
$N$ vertices, by the above reasoning $Q^{\circ}$ can be represented as a projection of a polytope 
with at most $N$ facets. Dualizing again and using that $\left(Q^{\circ}\right)^{\circ}=Q$, we 
conclude that $Q$ can be represented as a section of a polytope  with at most $N$ vertices. 
\bigskip
Approximations by projections are well suited for taking intersections and direct products.
\subhead (4.2) Operations on projections \endsubhead

For $i=1,2$, let $P_i \subset W_i$ be polyhedra, let $T_i: W_i \longrightarrow V$ be linear 
transformations, and let $Q_i=T_i\left(P_i\right)$.
Let $Q=\left(Q_1 \cap Q_2\right) \subset V$. We observe that $Q$ is the image of the polyhedron 
$P \subset W_1 \oplus W_2$, 
$$ P=\Bigl\{\left(x_1, x_2\right): \quad x_1 \in P_1,\quad  x_2 \in P_2, \quad \text{and} \quad
 T_1(x_1)=T_2(x_2) \Bigr\} $$
 under the linear transformation defined by $T\left(x_1, x_2\right)=T_1(x_1)$.
 If $P_i$ has at most $N_i$ facets for $i=1,2$ then $P$ has at most $N=N_1 + N_2$ facets.
 
For $i=1,2$, let $P_i \subset W_i$ be polyhedra, let $T_i: W_i \longrightarrow V_i$ be linear 
transformations, and let $Q_i=T_i\left(P_i\right)$.
Let $Q=\left(Q_1 \times Q_2\right) \subset V_1 \oplus V_2$. Then $Q$ is the image of the polyhedron
$P=\left(P_1 \times P_2\right) \subset W_1 \oplus W_2$ under the linear transformation defined
by $T\left(x_1, x_2\right)=\left(T_1\left(x_1\right), T_2\left(x_2\right)\right)$. Similarly, if 
$P_i$ has at most $N_i$ facets for $i=1,2$ then $P$ has at most $N=N_1+N_2$ facets.
\bigskip
One interesting feature of approximations by sections and projections is that it breaks symmetry.
Namely, to obtain a sufficiently close approximation of a {\it symmetric} convex body $B$ 
by a projection of a polytope $P$ with not too many facets, we may have to 
choose the polytope $P$ to 
be {\it non-symmetric}.

Let $O_d \subset {\Bbb R}^d$ be the standard octahedron (cross-polytope), see Section 2.3.
Clearly, $O_d$ has $2d$ vertices and hence there is a simplex $\Delta_{2d-1} \subset {\Bbb R}^{2d}$
with $2d$ vertices and a linear transformation $T: {\Bbb R}^{2d} \longrightarrow {\Bbb R}^d$ such 
that $T\left(\Delta_{2d-1}\right)=O_d$. In other words, $O_d$ can be represented as the 
projection of a polytope with $2d$ facets.

Suppose, however, that we want to construct a {\it symmetric} polytope $P \subset W$ and a 
linear transformation $T: W \longrightarrow {\Bbb R}^d$ such that $Q=T(P)$ approximates
$O_d$ within a factor of 2.
One can show that the number $N$ of facets of $P$ has to be exponentially large in $d$.

\proclaim{(4.3) Theorem} Let $O_d \subset {\Bbb R}^d$ be the cross-polytope 
$$O_d=\Bigl\{\left(\xi_1, \ldots, \xi_d \right): \ \sum_{i=1}^d |\xi_i| \leq 1 \Bigr\}$$
and suppose that $P \subset W$ is a symmetric polytope with $N$ facets and
$T: W \longrightarrow {\Bbb R}^d$ is a linear transformation such that 
$$Q \subset O_d \subset 2 Q \quad \text{for} \quad Q=T(P).$$ 
Then 
$$N \geq e^{cd} \quad \text{for some absolute constant} \quad c>0.$$
\endproclaim

\demo{Proof}
The proof uses the notion of the type 2 constant of a Banach space, see \cite{Pi89} and 
\cite{To89}.  

Let $V$ be a finite-dimensional vector space, let $B \subset V$ be a symmetric convex body 
and let $\| \cdot \|_B$ be the corresponding norm. The {\it type 2 constant} of $B$ is the 
smallest number $\kappa=\kappa(B)>0$ such that for  
any set of vectors $x_1, \ldots, x_m \in V$ we have 
$$\EE \Big\| \sum_{i=1}^m \epsilon_i x_i \Big\|_B^2 \leq \kappa^2(B) \sum_{i=1}^m \|x_i\|_B^2,$$
where the expectation is taken with respect to independent random signs $\epsilon_i$:
$$\epsilon_i=\cases \phantom{-}1 &\text{with probability} \ {1 \over 2} \\ -1 &\text{with probability} \ {1 \over 2}. 
\endcases$$
One can observe that $\kappa(B)=1$ if $B$ is the unit ball in some Euclidean metric
and that if 
$T: V \longrightarrow W$ is an invertible linear transformation then 
$\kappa(B)=\kappa\left(T(B)\right)$ for any symmetric convex body $B \subset V$.
 Since every convex body can be 
approximated by an ellipsoid, it follows that $\kappa(B)$ exists (in particular, is finite) for all symmetric 
convex bodies $B$.

Let $B \subset W$ be a symmetric convex body and let $U \subset W$ be a subspace. 
It is immediate that 
$$\kappa\left(B \cap U\right) \leq \kappa(B). \tag4.3.1$$

Furthermore, suppose that $T: W \longrightarrow V$ is a linear surjection. Then 
$$\kappa\bigl(T(B)\bigr) \leq \kappa(B). \tag4.3.2$$
To prove (4.3.2), we notice that 
$$\| Ty\|_{T(B)} \leq \|y\|_B \quad \text{for all} \quad y \in W$$
and that 
$$\text{for every \ } x\in V \ \text{there is \ } y \in W \ \text{such that} \ 
T(y)=x \ \text{and} \  \|y\|_B=\|x\|_{T(B)}.$$
To establish (4.3.2), let us choose any $x_1, \ldots, x_m \in V$ and then $y_1, \ldots, y_m \in W$ 
such that 
$$T(y_i)=x_i \quad \text{and} \quad \|y_i\|_B = \|x_i\|_{T(B)} \quad \text{for} \quad i=1, \ldots, m.$$
Then 
$$\EE \Big\|\sum_{i=1}^m \epsilon_i x_i \Big\|^2_{T(B)} \leq 
\EE \Big\| \sum_{i=1}^m \epsilon_i y_i\Big\|^2_B \leq \kappa^2(B) \sum_{i=1} \|y_i\|^2_B =
\kappa^2(B) \sum_{i=1}^m \|x_i\|^2_{T(B)}$$
and the proof of (4.3.2) follows.

Suppose now that $P$ is a symmetric polytope with $N$ facets. Then, up to a linear transformation,
$P$ can be represented as a section of an $N$-dimensional cube $I_N$, see Section 2.3.
Hence $\kappa(P) \leq \kappa\left(I_N\right)$ and if $Q$ is a projection of $P$, we have 
$$\kappa(Q) \leq \kappa(P) \leq \kappa\left(I_N\right).$$
Furthermore, under the conditions of the theorem,  we have 
$$\kappa\left(O_d\right) \leq 2\kappa(Q)  \leq 2\kappa\left(I_N\right).$$
On the other hand, one can show that 
$$\kappa\left(O_d\right) \geq \sqrt{d} \quad \text{and} \quad \kappa(I_N) \leq c \sqrt{\ln N}$$
for some absolute constant $c>0$, see Section 4 of \cite{To89}. This proves that $N$ has to be exponentially 
large in $d$.
{\hfill \hfill \hfill} \qed
\enddemo

Theorem 3.2 shows that the $d$-dimensional Euclidean ball is not approximated very well by a polytope with a subexponential in $d$ number of vertices or facets. The situation changes radically if 
we consider approximations by a projection of a polytope with a small number of facets or by 
a section of a polytope with a small number of vertices. For example, the intersection of 
a $2d$-dimensional octahedron $O_{2d}$ with a random $d$-dimensional subspace of ${\Bbb R}^{2d}$
approximates the $d$-dimensional Euclidean ball within an absolute constant, see, 
for example, \cite{Ba97}.
The results of \cite{B+89} on approximation of zonoids by zonotopes
 imply that for any $d$ and $\epsilon>0$ there is an  
$N=O^{\ast} \left(d \epsilon^{-2} \right) $, where $^{\ast}$ stands for some logarithmic factors, and a
linear
transformation $T: {\Bbb R}^N \longrightarrow {\Bbb R}^d$, such that the image $T\left(I_N\right)$
of the cube $I_N \subset {\Bbb R}^N$ approximates the unit ball $B \subset {\Bbb R}^d$ within 
a factor of $(1+\epsilon)$. Note that in this case we approximate a symmetric convex body by
the projection of a symmetric polytope.

A. Ben-Tal and A. Nemirovski \cite{BN01} provide a remarkable construction of a polytope
$P$ with $N=O\left(d \ln \epsilon^{-1} \right)$ of facets 
whose projection $T(P)$ approximates Euclidean ball $B \subset {\Bbb R}^d$ within a factor
of $(1+\epsilon)$. We sketch the construction below. Unlike in the case of zonotopal approximation
of \cite{B+89}, the polytope $P$ is not symmetric.

\subhead (4.4) A tight approximation of the Euclidean ball by a projection \endsubhead
Let 
$$B_d=\Bigl\{\left(\xi_1, \ldots, \xi_d \right): \quad \sum_{i=1}^d \xi_i^2 \leq 1 \Bigr\}$$
be the $d$-dimensional unit ball. The first idea of \cite{BN01} is to reduce the general case to that 
of $d=2$. Let us consider the $(d+1)$-dimensional round cone 
$$C_d=\Bigl\{\left(\xi_1, \ldots, \xi_d; \tau \right): \quad \sum_{i=1}^d \xi_i^2 \leq \tau^2 
\quad \text{and} \quad \tau \geq 0\Bigr\}.$$
One can observe that a close approximation of $C_d$ by a projection of a polyhedron with 
at most $N$ facets results in a close approximation of $B_d$ by a projection of a polytope 
with at most $N$ facets and vice versa.

Let us write $d=r+b$ for some positive integers $r$ and $b$ and let us consider the corresponding 
round cones of lower dimensions:
$$\split &C_r=\Bigl\{ \left(\xi_1, \ldots, \xi_r; \rho \right): \quad \sum_{i=1}^r \xi_i^2 \leq \rho^2 
\quad \text{and} \quad \rho \geq 0 \Bigr\} \quad \text{and} \\
&C_b=\Bigl\{ \left(\eta_1, \ldots, \eta_b ;\beta \right): \quad \sum_{i=1}^r \eta_i^2 \leq \beta^2 
\quad \text{and} \quad \beta \geq 0 \Bigr\}. \endsplit$$
In the space ${\Bbb R}^{d+3}$ with the coordinates $\xi_1, \ldots, \xi_b; \eta_1, \ldots, \eta_r; 
\rho, \beta, \tau$ we consider the subset
$$\split X=\Bigl\{&\left(\xi_1, \ldots, \xi_r; \eta_1, \ldots \eta_b; \rho, \beta, \tau\right): \\
&\sum_{i=1}^r \xi_i^2 \leq \rho^2, \quad \sum_{i=1}^b \eta_i^2 \leq \beta^2, \quad \rho^2 +\beta^2 
\leq \tau^2 \quad \text{and} \\
&\rho, \beta, \tau \geq 0 \Bigr\}. \endsplit$$
Clearly, $C_d$ is the image of $X$ under the projection which forgets $\rho$ and $\beta$.
On the other hand, we have 
$$X=\left(C_r \times C_b \times V_{\tau} \right) \cap \left(\overline{C_2} \times V_{\xi, \eta} \right),$$
where 
$$\overline{C_2} =\Bigl\{ \left(\rho, \beta, \tau\right): \quad \beta^2 + \rho^2 \leq \tau^2 \quad 
\text{and} \quad\rho, \beta, \tau \geq 0 \Bigr\}$$
is the ``quarter'' of the 3-dimensional round cone and  $V_{\tau}$ and $V_{\xi, \eta}$ are the appropriate coordinate subspaces.

Using the operations of direct product and intersection, see Section 4.2, one can show that 
good approximations of $C_r$, $C_b$, and $C_2$ by projections of polyhedra leads to a 
good approximation of $C_d$ by a projection of a polyhedron.

Essentially, the problem boils down to an efficient approximation of  the quarter of the disc:
$$\overline{B_2}=\Bigl\{\left(\xi, \eta \right): \quad \xi^2 + \eta^2 \leq 1 \quad \text{and} \quad 
\xi, \eta \geq 0 \Bigr\}.$$ 
For that, A. Ben-Tal and A. Nemirovski provide the following ingenious construction. 

Let us define the sequence of transformations $R_n$ of ${\Bbb R}^2$ by 
$$\split &\xi'=\xi \cos {\pi \over 2^n}  + \eta  \sin {\pi \over 2^n}  \\
               &\eta'=\Big|-\xi \sin {\pi \over 2^n} +\eta \cos {\pi \over 2^n} \Big|.
  \endsplit$$
Geometrically, $R_n$ is a clockwise rotation through an angle of $\pi/2^n$ followed by the reflection
in the $\xi$-axis if the obtained point lies in the lower halfplane. 
If we pick a point $z=(\xi, \eta)$ with $\xi, \eta \geq 0$ and apply the sequence of transformations
$R_2, R_3, R_4, \ldots, $ then the resulting sequence of points has its limit on the interval 
$0 \leq \xi \leq 1$, $\eta=0$ if and only if $z \in \overline{B_2}$.

Let us choose a positive integer $m$. In the space of variables $\xi_{k}, \eta_{k}$ for 
$k=1, \ldots, m$ we define the polyhedron $P_m$ 
by the equations and inequalities
$$\split &\xi_k=\xi_{k-1} \cos {\pi \over 2^k} + \eta_{k-1} \sin {\pi \over 2^k} \\
&\eta_k \geq -\xi_{k-1} \sin {\pi \over 2^k} +\eta_{k-1} \cos {\pi \over 2^k} \\
&\eta_k \geq \xi_{k-1} \sin {\pi \over 2^k} -\eta_{k-1} \cos {\pi \over 2^k} \quad \text{for} \quad 
k=2, \ldots, m\endsplit$$
and 
$$\xi_1, \eta_1 \geq 0, \quad 0 \leq \xi_m \leq 1, \quad \eta_m \leq {\pi \over 2^m}.$$
The projection of $P_m$ onto the coordinates $(\xi_1, \eta_1)$ approximates $\overline{B_2}$ within 
an error exponentially small in $m$.

\subhead (4.5) A general approximation construction \endsubhead 
Let $B \subset {\Bbb R}^d$ be a convex body we want to approximate by a projection of a polytope 
with at most 
$N$ facets, or, equivalently, by a section of a polytope with at most $N$ vertices. 
Here is a general construction. 

Without loss of generality, we assume that $B$ contains the origin in its interior and hence we
may view $B$ as the polar of its own polar $B^{\circ} \subset {\Bbb R}^d$. 
Next, we approximate $B^{\circ}$ by a sufficiently dense finite subset $X \subset B^{\circ}$ so 
that $X^{\circ}$ approximates $B$ well enough. 

Let us consider the space ${\Bbb R}^{X}$ of all functions $f: X \longrightarrow {\Bbb R}$ and 
let ${\Cal A} \subset {\Bbb R}^X$ be the affine subspace of all affine functions $f: X \longrightarrow {\Bbb R}$ whose average value on $X$ is 1. In other words, ${\Cal A}$ consists 
of the functions of the type $f(x)=\langle c, x \rangle + \alpha$ for some $c \in {\Bbb R}^d$ 
and $\alpha \in {\Bbb R}$ such that, additionally,
$${1 \over |X|} \sum_{x \in X} f(x)=1.$$ 
It is not hard to see that $X^{\circ}$ can be viewed as 
$$X^{\circ}=\Bigl\{f \in {\Cal A}: \quad f(x) \geq 0 \quad \text{for all} \quad x \in X \Bigr\}.$$
For a subset $Y \subset X$ let $\delta_Y \in {\Bbb R}^X$ be the indicator of $Y$, that is  
$$\delta_Y(x)=\cases 1 & \text{if \ } x \in Y \\ 0 &\text{if \ } x \notin Y.\endcases$$
We write simply $\delta_x$ instead of $\delta_{\{x\}}$.
Now we can write  
$$X^{\circ} ={\Cal A} \cap \co\left(\delta_x: \quad x \in X \right),$$
where ``$\co$'' stands for the {\it conic hull} of the set, that is, the set of all non-negative linear combinations of the elements of the set.

Let us choose a family ${\Cal F} \subset 2^X$. We define the approximation of $B_{\Cal F}$ as 
$$B_{\Cal F}={\Cal A} \cap \co\left(\delta_F: \quad F \in {\Cal F} \right).$$
Clearly, $B_{\Cal F}$ is the intersection of an affine subspace with a polyhedron with at 
most $|{\Cal F}|$ vertices. One can expect that the finer ${\Cal F}$ gets, the better approximation 
$B_{\Cal F}$ to $B$ we obtain. For example, we may choose $X$ as the intersection of $B$ with a 
sufficiently dense grid, $X=\left(B \cap \left(\epsilon {\Bbb Z}\right)^d \right)$ and then choose 
${\Cal F_k}$ consisting of the sets of points lying on an affine coordinate subspace of codimension 
$k$. Generally, we will have $|{\Cal F}_k|=O\left(d^k \epsilon^{-k} \right)$.

Using this approach, the second author was able to obtain the following approximation result 
for the Traveling Salesman Polytope $TSP_n$, see Example 1.2.

\proclaim{(4.6) Theorem} Let us choose an $\epsilon>0$ and let $V_n$ be the 
ambient space of $TSP_n$. Then there exists a 
polytope $P_n \subset W_n$ with at most $N=O\left(n^{4/\epsilon}\right)$ facets and
a linear transformation $T: W_n \longrightarrow V_n$ such that for $Q_n=T(P_n)$ 
one has
$$Q_n \subset TSP_n \subset (\epsilon n) Q_n.$$ 
\endproclaim

The proof of Theorem 4.6 is rather technical and will be presented elsewhere.

\subhead (4.7) Remarks and open questions \endsubhead 

Given a series of combinatorially defined polytopes $Q_n$, such as the Traveling Salesman Polytope $TSP_n$,
to construct a simpler polytope $P_n$ whose projection approximates $Q_n$, by now is a well-established
technique of ``lift and project'' in combinatorial optimization, see, for example, \cite{Ba01}. 
In particular, a remarkable success was achieved in showing that certain exponentially 
large (in $n$) families of facets of certain series of polytopes $Q_n$ can be obtained as 
projections of only polynomial in $n$ families of facets of $P_n$, see \cite{LS91}.
However, not much is known about how well the projections approximate {\it metrically}, even though 
various techniques were compared {\it combinatorially}, see \cite{La03}.

We also note that the authors of \cite{BN01}, while constructing their remarkable approximation 
of the Euclidean ball by the projection of a polytope with a small number of facets (cf. Section 4.4), were 
motivated by very practical questions. Namely, they used their approximation to reduce 
convex quadratic programming problems to linear programming problems (see \cite{BN01} for 
details) and hence use the linear programming solver available to them to solve quadratic 
programs. 

Some basic questions regarding approximations by projections remain unanswered.
\smallskip
(4.7.1) {\it Obstructions to being a projection.} Let $Q \subset V$ be a polytope. How can one
possibly prove that $Q$ {\it cannot be} the projection of a polytope $P$ with at most $N$ facets?
In the case when $Q$ is symmetric and $P$ is required to be symmetric, a possible argument 
goes via the type 2 constant, see Theorem 4.3. Another example is the proof of \cite{Ya91}
that $TSP_n$ cannot be a projection of a polytope $P_n$ with a polynomial in $n$ number 
of facets provided the projection respects the symmetries of $TSP_n$. In general, if 
the lifting $P$ is not required to be symmetric, no viable argument seems to be known. It would
be interesting to find out whether some appropriate notion of a non-symmetric type can be of help.
\smallskip
(4.7.2) {\it The quality of approximation of a general body.} How well can a general $d$-dimensional convex body $B$ be approximated by 
the projection of a polytope with at most $N$ facets? In particular, the following question seems 
to be of interest. Suppose that $B$ is symmetric and that the projection of a polytope $P$ approximates
$B$ within a factor of 2. Is it true that in the worst case the number $N$ of facets of $P$ should be at least 
exponential in $d$, $N \geq c^d$ for some absolute constant $c>1$? In Section 5 we discuss a
certain construction which suggests that maybe the exponential bound can be broken and that 
we can have $N =O^{\ast}\left(c^{\sqrt{d}}\right)$, where, as usual, $\ast$ stands for some 
logarithmic factors.
\smallskip
(4.7.3) {\it Approximation of the $l^p$ ball.} Let $p \geq 1$ and let
$$B(d,p)=\Bigl\{ \left(\xi_1, \ldots, \xi_d\right): \quad \sum_{i=1}^d \xi_i^p \leq 1 \Bigr\}$$
be the unit ball in the $l^p$ norm. One can observe that the first step of the construction of 
 Section 4.4 can be extended to $B(d,p)$ thus providing an approximation
of $B(d,p)$ within a factor of $(1+\epsilon)$ by the projection of a polytope with at most
$N=O\left(d \epsilon^{-1}\right)$ facets. It is not clear whether for a general $p$ one 
can replace $\epsilon^{-1}$ by $\ln \epsilon^{-1}$, though this is definitely the case for $p=1,2$,
and $+\infty$.

\head 5. A ``soft'' approximation of a symmetric convex body \endhead 

Let  $B \subset V$ be a symmetric convex body. We identify $V^{\ast}$ with the subspace of linear 
functions in the space $C\left(B\right)$ of all {\it continuous} functions $f: B \longrightarrow {\Bbb R}$.
Then $B^{\circ}$ can be identified with the set of  {\it linear} functions 
$f: B \longrightarrow {\Bbb R}$ such that $f(x) \leq 1$ for all $x \in B$. 

In this section, we prove the following main result.
\proclaim{(5.1) Theorem} 
 Let $B \subset V$ be a symmetric $d$-dimensional convex body  and let 
$C\left(B\right)$ be the space of all continuous functions on $B$.
Then there exists a polytope $R \subset C\left(B\right)$ such that the following holds.
\roster
\item For all $h \in R$ we have 
$$h(x ) \leq 1 \quad \text{for all} \quad  x  \in B;$$
\item For any $\ell \in B^{\circ}$ there exists a function 
$h \in R$ such that 
$$\big| \ell(x) -h(x) \big| \leq  \gamma \ell^2(x)  \quad \text{for all} \quad x \in B$$
and some absolute constant $\gamma$.
\item The polytope $R$ has at most $\exp\left\{ \alpha\sqrt{d} \ln d \right\}$ vertices
for some absolute constant $\alpha$. 
\endroster
\endproclaim

If we could claim that the function $h$ in Part (2) of the theorem is {\it linear}, then we must have 
had $h(x)=\ell(x)$ in Part (2) and we would have obtained the representation 
$$B^{\circ}=R \cap V^{\ast}$$
of $B^{\circ}$ as the intersection of a polytope with at most $\exp\left\{\alpha \sqrt{d} \ln d \right\}$ 
vertices and a subspace $V^{\ast}$. By duality, that would have implied that $B$ is the projection of 
a polytope with at most $\exp\left\{\alpha \sqrt{d} \ln d \right\}$ facets. In general, however, $h$ is not a linear function, but, as will 
follow from the proof, is a piecewise polynomial. If $\ell \in \epsilon B^{\circ}$ for some $0< \epsilon <1$ 
then $h$ approximates $\ell$ within an error of $O\left(\epsilon^2 \right)$, so the points of $B^{\circ}$ that 
are closer to the origin are better approximated. Intuitively, we obtain a set close to $B^{\circ}$ if 
we slightly ``bend'' $V^{\ast}$ and then intersect it with $R$.
 
\demo{Proof of Theorem 5.1} Since $B$ is symmetric, we can find an ellipsoid $E \supset B$ centered at the origin 
such that 
$${1 \over \sqrt{d}}E \subset B \subset E.$$ 
On the other hand, approximating the ellipsoid by the projection of a polytope (see Section 4.4), one 
can construct a polytope $P \subset W$ with $N =O(d)$ facets and a linear transformation $T: W \longrightarrow V$ such that 
$${1 \over 2} T(P) \subset E \subset T(P).$$  
Summarizing,
$${1 \over 2 \sqrt{d}} T(P) \subset B \subset T(P). \tag5.1.1$$ 
Without loss of generality, we assume that $P \subset W$ is full-dimensional 
and contains the origin in its interior. 
Suppose that 
$$P=\Bigl\{w \in W: \quad g_i(w) \leq 1, \ i=1, \ldots, N \Bigr\}, \quad 
\text{where} \quad g_i: W \longrightarrow {\Bbb R} $$
are linear functions. Let us choose the smallest positive integer $k > 2\sqrt{d}$. For 
a multiset $I$ of numbers from 1 to $N$ of cardinality at most $k$ (counting multiplicities), we let 
$$g_I=1-\prod_{i \in I} \left(1-g_i\right).$$
Hence $g_I: W \longrightarrow {\Bbb R}$ are polynomials, $\deg g_I \leq k$.
It is immediate that $g_I(w) \leq 1$ for all $w \in P$. 

Suppose now that $\ell \in B^{\circ}$, so $\ell: B \longrightarrow {\Bbb R}$ is a linear function
such that $\ell(x) \leq 1$ for all $x \in B$. In view of (5.1.1), we have $\ell(x) \leq k$ for all 
$x \in T(P)$.
Let $f: W \longrightarrow {\Bbb R}$ be the lifting of $\ell$ defined by $f(w)=
\ell\left(T(w)\right)$ for $w \in W$. Hence $f: W \longrightarrow {\Bbb R}$ is a linear
function and $f(w) \leq k$ for all $w \in P$. Therefore, $k^{-1} f(w) \leq 1$ 
for all $w \in P$
and hence 
$$k^{-1} f \in \conv\Bigl(0, \ g_i: \quad i=1, \ldots, N \Bigr).$$ 
Therefore,
$$F=1-\left(1-k^{-1} f\right)^k \in \conv\Bigl(0, \ g_I: \quad  |I| \leq k \Bigr). \tag5.1.2$$
Now, since $B$ symmetric, we have $|\ell(x)| \leq 1$ for all $x \in B$. Therefore, 
for all $w \in P$ such that $T(w) \in B$, we have $|f(w)| \leq 1$ and hence
$$|F(w) - f(w)| \leq \gamma f^2(w) \quad \text{for all} \quad w \in P \quad 
\text{such that} \quad T(w) \in B. \tag5.1.3$$

Now we are ready to define $R \subset C(B)$. Let us fix a scalar product in $W$ and hence 
the Lebesgue measure on every affine subspace of $W$. Let us define $h_I \in C(B)$ by 
$$h_I(x)=\text{the average value of \ } g_I(w) \quad \text{over all} \quad w \in P \cap T^{-1}(x).$$
 Let us define $R$ to be the convex hull of the origin and all 
the functions $h_I$ as $I$ ranges over all multisets with the elements from $\{1, \ldots, N\}$ 
and of cardinality at most $k$, counting multiplicities. Clearly, the number of vertices of $R$ does not exceed
$N^k$, so Part(3) follows. 

Since $h_I(x) \leq 1$ for all $I$ and all $x \in B$, Part (1) follows as well. 

To prove Part (2) we choose
$$h(x)=\text{the average value of \ } F(w) \quad \text{over all} \quad w \in P \cap T^{-1}(x),$$
where $F$ is defined by (5.1.2). Clearly, $h \in R$ and Part (2) follows by (5.1.3).
{\hfill \hfill \hfill} \qed
\enddemo
\bigskip
Let $\ell \in V^{\ast}$ be a linear function $\ell: B \longrightarrow {\Bbb R}$. There seems to be no 
efficient way to check whether there is a function $h \in R$ such that $|\ell(x)-h(x)| \leq \gamma \ell^2(x)$
for all $x \in B$. We can relax the condition by replacing the uniform distance by the distance 
in the $L^2(B, \mu)$ norm for some Borel probability measure $\mu$ on $B$:
$$\|f\|_2=\left(\int_B f^2 \ d \mu\right)^{1/2}.$$
Then, checking whether for a given $\ell \in V^{\ast}$ there is an $h \in R$ 
such that 
$$\|\ell-h\|_2 \leq \gamma \|\ell^2\|_2$$  becomes a problem 
of convex quadratic programming which can be solved  roughly in 
$\exp\left\{O\bigl(\sqrt{d} \ln d\bigr)\right\}$ time (in the bit model).

Let us choose a sufficiently small $\epsilon>0$. Let us ``accept'' a given linear function
$\ell \in V^{\ast}$ if there exists an $h \in R$ such that 
$$\|\ell-h\|_2 \leq \gamma \epsilon \|\ell\|_2$$
and ``reject'' it otherwise. Hence we get an algorithm of
$\exp\left\{O\bigl(\sqrt{d} \ln d\bigr)\right\}$ complexity such that given a linear function 
$\ell:B \longrightarrow {\Bbb R}$
\bigskip
(i) the algorithm accepts $\ell$ if $\ell(x) \leq \epsilon$ for all $x \in B$;
\smallskip
(ii) if the algorithm accepts $\ell$ then there is a function $h: B \longrightarrow {\Bbb R}$ with 
$h(x) \leq 1$ for all $x \in B$ and $\|\ell-h\|_2 \leq \epsilon \gamma \|\ell\|_2 $, where $\gamma$ is 
an absolute constant.
\bigskip 
This is somewhat similar to the situation of the ``property testing'' in computer science, see, for example, 
\cite{Go99}.  

One can observe that the estimates of Theorem 5.1 can be extended 
in a number of ways. If we know that $B$ can be approximated by an ellipsoid within a factor 
$\rho \geq 1$ then we can construct a polytope $R$ with the number of vertices not exceeding
$\exp\left\{ \alpha \rho \ln d \right\}$. More generally, we can require that $R$ has not more than
$\exp\left\{\alpha k \ln d\right\}$ vertices for $k \leq \rho$
 if we replace the condition (2) in Theorem 5.1 by 
$|k \rho^{-1} \ell(x)-h(x)| \leq \gamma k^2 \rho^{-2} \ell^2(x)$.

\head 6. Approximations by a section of the cone of positive semidefinite quadratic forms 
\endhead

Let $W$ be a finite-dimensional real vector space and let $S(W)$ be the space of all 
quadratic forms $q: W \longrightarrow {\Bbb R}$. Let $S_+(W) \subset S(W)$ be the 
convex cone of positive semidefinite quadratic forms, that is, the quadratic forms such that 
$q(w) \geq 0$ for all $w \in W$. The membership problem for the cone $S_+(W)$ can be 
solved in time polynomial in $\dim W$ both in the algebraic and the bit models of computation.
Roughly, checking that $q$ is positive semidefinite can be done by a certain variation 
of the Sylvester criterion. Identifying $W={\Bbb R}^d$, we conclude that $q$ is positive 
semidefinite if and only if the form
$$q_{\epsilon}(x)=q(x) + \epsilon \left(\xi_1^2 + \ldots + \xi_d^2 \right)$$
is strictly positive definite for all $\epsilon>0$. Sylvester's criterion then implies that the 
the $d$ principle minors of $q_{\epsilon}$ should be positive for all $\epsilon>0$, which, 
in turn, reduces to checking that $d$ univariate polynomials in $\epsilon$ are positive 
for all $\epsilon>0$.    

Let $B \subset V$ be a $d$-dimensional convex body. We may try to approximate $B$ by a set $X$ which is 
the intersection of the cone $S_+(W)$ for some $W$ and a $d$-dimensional affine subspace 
identified with $V$. The main accomplishments of this approach are associated with 
approximations of the cut polytope, see \cite{DL97}.

\subhead (6.1) Approximating the cut polytope \endsubhead 
For $n$-vectors $x=\left(\xi_1, \ldots, \xi_n\right)$ and 
$y=\left(\eta_1, \ldots, \eta_n \right)$, let us define the $n \times n$ matrix $x \otimes y$ as the 
matrix with the $(i,j)$th entry equal to $\xi_i \eta_j$ and let $\Mat_n$ be the vector 
space of $n \times n$ matrices. 

Let us identify the space $S\left({\Bbb R}^n\right)$ of quadratic forms on ${\Bbb R}^n$ with the 
space of $n \times n$ symmetric matrices. 
The {\it cut polytope} $CUT_n \subset S\left({\Bbb R}^n\right)$ is defined as the convex hull of all $n \times n$ matrices 
$x \otimes x$ for all vectors $x=\left(\xi_1, \ldots, \xi_n \right)$ with $\xi_i=\pm1$.
As is the case with the Traveling Salesman Polytope, the membership problem for the cut 
polytope is NP-complete, cf. \cite{DL97}.

We consider $CUT_n$ as a subset of the space $S\left({\Bbb R}^n \right)$ of 
symmetric $n \times n$ matrices. Let ${\Cal A} \subset S\left({\Bbb R}^n\right)$ be  the affine subspace 
consisting of the matrices with 1's on the diagonal. It turns out that  the intersection 
$S_+\left({\Bbb R}^n \right) \cap {\Cal A}$ approximates $CUT_n$ within a logarithmic factor with 
respect to the center at the identity matrix $I \in S\left({\Bbb R}^n \right)$:
$$CUT_n \  \subset \  {\Cal A} \cap S_+\left({\Bbb R}^n\right)\  \subset \ c \ln n \left( CUT_n\right)$$  
for some absolute constant $c>0$. The logarithmic factor cannot be improved, see 
\cite{A+06}.

A variation of $CUT_n$ is what we call the {\it asymmetric cut polytope} $ACUT_n$ defined 
as the convex hull of all matrices $x \otimes y$ where $x$
and $y$ are $n$-vectors with the coordinates $\xi_i, \eta_j =\pm 1$.  
Again, the membership problem for $ACUT_n$ is  NP-complete. It turns out that 
$ACUT_n$ can be tightly approximated by the {\it projection} of a section of a cone 
of positive semidefinite quadratic forms. Namely, let ${\Cal A} \subset S\left({\Bbb R}^{2n}\right)$ 
be the affine subspace of symmetric $2n \times 2n$ matrices with 1's on the diagonal. 
Let $\phi: S\left({\Bbb R}^{2n} \right) \longrightarrow \Mat_n$ be the projection
$$\phi \left( \matrix B & X\\ X & C \endmatrix \right) = X$$
and let 
$$Q_n=\phi\left({\Cal A} \cap S\left({\Bbb R}^{2n} \right)\right).$$
In words: $Q_n$ is the set of all possible  $n \times n$ upper right corner
 submatrices of a $2n \times 2n$ 
 positive semidefinite matrix with 1's on the diagonal.
It is not hard to see that $ACUT_n \subset Q_n$. Indeed, the matrix $z \otimes z$, where 
$z=\left(\xi_1, \ldots, \xi_n, \eta_1, \ldots, \eta_n \right)$ with $\xi_i, \eta_j =\pm 1$ 
is positive semidefinite with 1's on the diagonal and has the matrix $x \otimes y$ for  
$x=\left(\xi_1, \ldots, \xi_n\right)$  and $y=\left(\eta_1, \ldots, \eta_n \right)$ as 
its upper right corner submatrix. 

It turns out that $Q_n$ approximates 
$ACUT_n$ within a constant factor:
$$Q_n\ \subset\  ACUT_n \  \subset \  \kappa Q_n$$
for some absolute constant $\kappa$, called the {\it Grothendieck constant}. Its exact value is 
not known, but it is known that
$$1.5708 \approx {\pi \over 2} \leq \kappa \leq {\pi \over 2 \ln \left(1+\sqrt{2}\right)} \approx 1.7822,$$
see \cite{A+06} and \cite{AN06} for survey and recent developments.

We note that the membership problem for the projection of a section of a cone of 
positive semidefinite quadratic forms is an instance of {\it semidefinite programming},
which can be solved in polynomial time (in the bit model), though only approximately, cf. \cite{Kl02}.
\bigskip
Now we describe a construction for approximating a general convex body by a section of the 
cone of positive semidefinite forms. 

Let $B \subset V$ be a convex body containing the origin in its interior. As in Section 2, we think of $B$ as the polar $B=\left(B^{\circ}\right)^{\circ}$ to its own polar $B^{\circ} \subset V^{\ast}$. Hence 
$$B=\Bigl\{x \in V: \quad \ell(x) \leq 1 \quad \text{for all} \quad \ell \in B^{\circ} \Bigr\}.$$
Let us choose a positive integer $k$ and let $W_k$ be the space of all polynomials 
$p: V^{\ast} \longrightarrow {\Bbb R}$, $\deg p \leq k$. In particular, for any fixed $k$, the dimension 
of $W_k$ is bounded by a polynomial in the dimension of $V$.

Let us choose a Borel probability measure $\mu$
on $B^{\circ}$. 
For a point $v \in V$, let $q_v: W_k \longrightarrow {\Bbb R}$ be the quadratic form defined by
$$q_v(p)=\int_{B^{\circ}} \bigl(1-\ell(v)\bigr) p^2(\ell) \ d \mu(\ell) \quad \text{for} 
\quad p \in W_k.$$
Clearly, if $v \in B$ then the form $q_v$ is positive definite and the set
$$\Bigl\{q_{v}: \  v \in V \Bigr\} \subset S\left(W_k\right)$$ 
is an affine subspace. This allows us to define an approximation $X_k$ to $B$ by 
$$X_k=\Bigl\{v \in V:\ q_v \in S_+\left(W_k\right) \Bigr\}.$$
Thus $B \subset X_k$ for all $k$. To show that $X_k$ approximates $B$ reasonably well,
we have to show that for any point $v \notin B$ which is sufficiently far away from $B$ we
can find a polynomial $p: V^{\ast} \longrightarrow {\Bbb R}$ which takes large values on 
$\ell \in B^{\circ}$ such that $\ell(v)>1$ and  small values everywhere else on $B^{\circ}$.
Besides, the value of $\mu\bigl\{\ell \in B^{\circ}: \quad \ell(v) > 1 \bigr\}$ should be sufficiently 
large, in particular, the exterior angle measure $\mu_{B^{\circ}}$ discussed in Section 2.4 can be 
of help. 

In the case of the Traveling Salesman Polytope, the second author obtained the following result.
\proclaim{(6.2) Theorem} Let us choose an $\epsilon>0$ and let $V_n$ be the ambient space 
of $TSP_n$. Let $TSP^{\circ}_n \subset V_n^{\ast}$ be the polar of the 
Traveling Salesman Polytope with respect to its center as the origin. Then, there exists a set 
$Q_n \subset V_n^{\ast}$ isometric to the section of the cone of positive semidefinite quadratic forms in
$n^{O(1/\epsilon)}$ variables by an affine subspace such that 
$$Q_n \subset TSP_n^{\circ} \subset (\epsilon n) Q_n.$$
\endproclaim

The proof is presented in \cite{Ve06}. It is not clear whether the bound is sharp.
For example, the vertices of $TSP_n^{\circ}$ corresponding to the standard facets $x_{ij}=0$ of the
$TSP$ actually belong to $Q_n$ while the vertices of $TSP_n^{\circ}$ corresponding to 
the {\it subtour elimination facets}, see Chapter 58 of \cite{Sc03}, lie in $\alpha_n Q_n$ for $\alpha_n=O\left(\sqrt{n}\right)$.

\head Acknowledgments \endhead

The authors are grateful to G. Schechtman for pointing out to connections between inapproximability
and the type constants of Banach spaces (Theorem 4.3) and to A. Nemirovski for explaining 
a way to tightly approximate the Euclidean ball by the projection of a polytope with not too
many facets (Section 4.4) and for encouragement.

\enddocument
\end